\documentclass{amsart}

\usepackage{cite}
\usepackage{textcomp}
\usepackage{amsmath, amssymb, amsfonts}
\usepackage{dsfont}
\usepackage{accents}
\usepackage{algorithm}
\usepackage{algorithmic}
\usepackage{color}
\usepackage{graphicx}
\usepackage{booktabs}
\usepackage{caption}
\usepackage{subcaption}
\usepackage{enumerate}
\usepackage{mathrsfs}
\usepackage{mathtools}

\DeclareMathOperator{\tr}{tr}

\DeclareMathOperator{\diag}{diag}
\DeclareMathOperator{\voi}{VoI}
\DeclareMathOperator{\E}{\mathsf{E}}
\DeclareMathOperator{\Cov}{\mathsf{cov}}
\DeclareMathOperator{\ProbM}{\mathsf{P}}
\DeclareMathOperator{\Prob}{\mathsf{p}}
\DeclareMathOperator{\argmin}{arg min}

\newtheorem{definition}{Definition}
\newtheorem{lemma}{Lemma}
\newtheorem{theorem}{Theorem}
\newtheorem{proposition}{Proposition}

\newtheorem{remark}{Remark}

\begin{document}
\title{Value of Information in Feedback Control: Quantification}
\author{
Touraj Soleymani, John S. Baras, and Sandra Hirche}
\thanks{Corresponding Author: Touraj Soleymani (touraj@kth.se). Journal: \emph{IEEE Transactions on Automatic Control}.}

\maketitle

\begin{abstract}
Although transmission of a data packet containing sensory information in a networked control system improves the quality of regulation, it has indeed a price from the communication perspective. It is, therefore, rational that such a data packet be transmitted only if it is valuable in the sense of a cost-benefit analysis. Yet, the fact is that little is known so far about this valuation of information and its connection with traditional event-triggered communication. In the present article, we study this intrinsic property of networked control systems by formulating a rate-regulation tradeoff between the packet rate and the regulation cost with an event trigger and a controller as two distributed decision makers, and show that the valuation of information is conceivable and quantifiable grounded on this tradeoff. In particular, we characterize an equilibrium in the rate-regulation tradeoff, and quantify the value of information ${\voi}_{{\mathit{k}}}$ there as the variation in a so-called value function with respect to a piece of sensory information that can be communicated to the controller at each time~${\mathit{k}}$. We prove that, for a multi-dimensional Gauss--Markov process, ${\voi}_{{\mathit{k}}}$ is a symmetric function of the discrepancy between the state estimates at the event trigger and the controller, and that a data packet containing sensory information at time ${\mathit{k}}$ should be transmitted to the controller only if ${\voi}_{{\mathit{k}}}$ is nonnegative. Moreover, we discuss that ${\voi}_{{\mathit{k}}}$ can be computed with arbitrary accuracy, and that it can be approximated by a closed-form quadratic function with a performance guarantee.

\end{abstract}

\section{Introduction}
Networked control systems are spatially distributed systems wherein feedback control loops are closed over communication channels~\cite{baillieul2007X}. Commonly, in a networked control system, data packets containing sensory information are transmitted to the controller in a periodic way as this facilitates the analysis of such a system~\cite{alur2007}. It has, however, been conceived that not every one of these data packets has the same effect on the system performance, and that one should employ a mechanism, i.e., event trigger, that transmits a data packet \emph{only when a significant deviation in the system occurs}~\cite{Astrom:2002eg}. This adaptive communication has a major consequence: a dramatic reduction in the number of packet transmissions guaranteeing some level of system performance, which has been found appealing, and led to extensive development of event-triggered systems in different contexts even beyond control including consensus~\cite{dimos2012}, fault detection~\cite{li2017}, optimization~\cite{meinel2014}, and signal processing~\cite{tsividis2010}.

Although transmission of a data packet containing sensory information in a networked control system decreases the uncertainty of the controller and improves the quality of regulation, it has indeed a price from the communication perspective. It is, therefore, rational that such a data packet be transmitted only if it is valuable in the sense of a cost-benefit analysis, i.e., \emph{only if its benefit surpasses its cost}. Yet, the fact is that little is known so far about this valuation of information and its connection with the above-mentioned adaptive communication. In the present article, we study this intrinsic property of networked control systems by formulating a \emph{rate-regulation tradeoff} between the packet rate and the regulation cost, and show that the valuation of information is conceivable and quantifiable grounded on this tradeoff.

The rate-regulation tradeoff in our study involves a stochastic optimization problem with an event trigger and a controller as two distributed decision makers. Unfortunately, this problem for the joint design of the event trigger and the controller is intractable (see e.g., \cite{wu2013, ramesh2013}). The reasons are that, in general, the underlying information structure is non-classical, the optimal estimator at the controller is nonlinear with no analytical solution due to a signaling effect, and estimation and control are coupled due to a dual effect. Nevertheless, in this article, we characterize an equilibrium at which neither decision maker has a unilateral incentive to change its policy, and quantify the \emph{value of information} $\voi_k$ there as the variation in a so-called value function with respect to a piece of sensory information that can be communicated to the controller at each time~$k$. We study the issue of global optimality of this very equilibrium in a separate article~\cite{voi2}.

We argue that the value of information systematically captures the \emph{semantics} of data packets by taking into account their potential impacts, and that a strategy based on the value of information optimally shapes the information flow in networked control systems. As such, the value of information can be regarded as a semantic metric that determines the \emph{right piece of information}, a concept that is not defined in classical data communication, while it is crucial to the development of future communication networks~\cite{uysal2021}. Note that previously Dempster~\cite{dempster1981} and Davis~\cite{davis1989, davis1991} studied the value of information in the context of optimal control. However, in these works, the value of information was defined as the variation in a value function with respect to relaxation of the non-anticipativity constraint at the controller. It is obvious that our perspective here is fundamentally different\footnote{For the preliminary work of the authors on the topic of the value of information in feedback control, see e.g., \cite{mywodespaper, soleymaninecsys, soleymani2016-cdc}.}.

\subsection{Related Work}
An event trigger can be used at the sensor side to reduce the number of packet transmissions in the observation channel\footnote{The observation channel is a communication channel that is placed between the sensor and the controller. In contrast, the command channel is a communication channel that is placed between the controller and the actuator.}, or at the controller side to reduce that in the command channel. We are here interested in finding the optimal decision policies in the former case\footnote{A control problem with data-rate constraints on both observation and command channels can in effect be converted to one with a data-rate constraint only on the observation channel (see~e.g.,~\cite{nair2000}).}, where the event trigger and the controller become two distributed decision makers. To elucidate the essence of the underlying problem, we suppose that network-induced effects such as quantization, packet dropouts, and time-varying delays are negligible. In this context, {\AA}str\"om and Bernhardsson~\cite{Astrom:2002eg} showed that for a scalar linear diffusion process, with impulse control and under a sampling rate constraint, event-triggered sampling outperforms periodic sampling in the sense that it attains a lower mean error variance. 

In addition, several works have addressed optimal event-triggered estimation, and found optimal triggering policies~\cite{imer2010, rabi2012, lipsa2011, molin2017}. Notably, Imer and Ba{\c{s}}ar \cite{imer2010} studied the optimal event-triggered estimation of a scalar Gauss--Markov process with perfect information\footnote{Perfect information refers to a situation where the exact value of the state of the process can be observed at each time. In contrast, imperfect information refers to a situation where only a noisy version of the output of the process can be observed at each time.} based on dynamic programming by assuming that the triggering policy is symmetric threshold, and obtained the optimal threshold value of the policy. Rabi~\emph{et~al.}~\cite{rabi2012} formulated the optimal event-triggered estimation of the scalar Wiener and scalar Ornstein--Uhlenbeck processes with perfect information as an optimal multiple stopping time problem, and showed that the optimal triggering policy is symmetric threshold when negative information\footnote{Negative information refers to any information that can be inferred associated with non-transmission by the receiver/controller.} is discarded. Lipsa and Martins~\cite{lipsa2011} used majorization theory to analyze the optimal event-triggered estimation of a scalar Gauss--Markov process with perfect information, and proved that the optimal triggering policy is symmetric threshold despite the presence of negative information. Later, Molin and Hirche~\cite{molin2017} studied the convergence properties of an iterative algorithm for the optimal event-triggered estimation of a scalar Markov process with perfect information and symmetric noise distribution, and found a result coinciding with that in~\cite{lipsa2011}.

In the joint design of the event trigger and the controller, a separation between estimation and control is not given a priori. Therefore, the above results on optimal event-triggered estimation do not apply directly to optimal event-triggered control. There exist, however, a number of works that have specifically addressed optimal event-triggered control, and found optimal control policies~\cite{molin2013, ramesh2013, demirel2018}. In particular, Molin and Hirche~\cite{molin2013} investigated the optimal event-triggered control of a Gauss--Markov process with perfect information, and showed that the optimal control policy is certainty equivalent when the triggering policy is reparametrizable in terms of primitive random variables. Ramesh~\emph{et~al.}~\cite{ramesh2013} studied the dual effect in the optimal event-triggered control of a Gauss--Markov process with perfect information, and proved that the dual effect generally exists. They also proved that the certainty equivalence principle holds if and only if the triggering policy is independent of the control policy. Later, Demirel~\emph{et~al.}~\cite{demirel2018} addressed the optimal event-triggered control of a Gauss--Markov process with imperfect information by adopting a stochastic triggering policy that preserves the Gaussianity of the conditional distribution, and showed that the optimal control policy remains certainty equivalent.

\subsection{Contributions and Outline}
In this article, we introduce the notion of the value of information, and establish a theoretical framework for its quantification. More specifically, we prove the existence of an equilibrium in the rate-regulation tradeoff for a multi-dimensional Gauss--Markov process with imperfect information without any limiting assumptions on the information structure or the policy structure, and quantify $\voi_k$ at this equilibrium, where the optimal estimator at the controller becomes linear, the design of the event trigger and the controller becomes separated, and the control becomes neutral. We prove that $\voi_k$ is a symmetric function of the discrepancy between the state estimates at the event trigger and the controller, and that a data packet containing sensory information should be transmitted to the controller at time $k$ only if $\voi_k$ is nonnegative. Moreover, we discuss that $\voi_k$ can be computed with arbitrary accuracy, and that it can be approximated by a closed-form quadratic function with a performance guarantee.

In our analysis, we show that a symmetric threshold triggering policy based on the value of information and a certainty-equivalent control policy based on a non-Gaussian linear estimator are in fact person-by-person optimal. This structural result applies to multi-dimensional Gauss--Markov processes. Therefore, it is in contrast with the results in \cite{lipsa2011,molin2017}, which are restricted to scalar Gauss--Markov processes. Our triggering policy in its specialized scalar form, however, is consistent with the one obtained in \cite{lipsa2011,molin2017}, and requires similar complexity for the computation of the threshold. In addition, the above structural result determines the triggering policy and the control policy jointly at an equilibrium, and asserts that the conditional mean used within the control policy, as we will see, is not affected by negative information at all. Hence, it is different from the results in~\cite{molin2013, ramesh2013}, which specify only the optimal control policy when the triggering policy is fixed, providing no insight into the associated conditional mean in the optimal design when the triggering policy is not fixed a priori.

The remainder of the article is organized in the following way. We formulate the rate-regulation problem in Section~\ref{c1:problem-formulation}, and present our results on the characterization and the computation of the value of information in Section~\ref{c1:main-results1}. We then provide our numerical examples in Section~\ref{c1:examples}. Finally, we make concluding remarks in Section~\ref{c1:conclusion}.

\subsection{Preliminaries}
In the sequel, the sets of real numbers and non-negative integers are denoted by $\mathbb{R}$ and $\mathbb{N}$, respectively. For $x,y \in \mathbb{N}$ and $x \leq y$, the set $\mathbb{N}_{[x,y]}$ denotes $\{z \in \mathbb{N} | x \leq z \leq y\}$. For matrices $X$ and $Y$, the relations $X \succ 0$ and $Y \succeq 0$ denote that $X$ and $Y$ are positive definite and positive semi-definite, respectively. The indicator function of a subset $\mathcal{A}$ of a set $\mathcal{X}$ is denoted by $\mathds{1}_\mathcal{A}:\mathcal{X} \to \{0,1\}$.  The probability measure of a random variable $x$ is represented by $\ProbM(x)$, its probability density or probability mass function by $\Prob(x)$, and its expected value and covariance by $\E[x]$ and $\Cov[x]$, respectively.

\begin{definition}[Dual effect]
For a given control system, let $\mathcal{I}_k^c$ be the information set of the controller at time $k$, and $\tilde{\mathcal{I}}_k^c$ be the information set of the controller at time $k$ when all control inputs are equal to zero. The control has no dual effect of order $r$, $r \geq 2$, (see e.g.,~\cite{bar1974dual}) if
\begin{align*}
	\E[M_{k,i}^r | \mathcal{I}_k^c ] = \E[M_{k,i}^r | \tilde{\mathcal{I}}_k^c],
\end{align*}
where $M_{k,i}^r = (x_{k,i} - \E[x_{k,i} | \mathcal{I}_k^c])^r$ is the $r$th central moment of the $i$th component of the state $x_k$ conditioned on $\mathcal{I}_k^c$. In other words, the control has no dual effect if the expected future uncertainty is not affected by the prior control inputs.
\end{definition}

\begin{definition}[Nash equilibrium]
For a given team game with two decision makers, let $\gamma^1 \in \mathcal{G}^1$ and $\gamma^2 \in \mathcal{G}^2$ be the decision policies of the decision makers, where $\mathcal{G}^1$ and $\mathcal{G}^2$ are the sets of admissible policies, and $L(\gamma^1,\gamma^2)$ be the associated loss function. A policy profile $(\gamma^{1\star},\gamma^{2\star})$ represents a Nash equilibrium (see e.g., \cite{yuksel}) if
\begin{align*}
	L(\gamma^{1\star},\gamma^{2\star}) \leq L(\gamma^{1},\gamma^{2\star}), \ \text{for all } \gamma^1 \in \mathcal{G}^1,\\[2.25\jot]
	L(\gamma^{1\star},\gamma^{2\star}) \leq L(\gamma^{1\star},\gamma^{2}), \ \text{for all } \gamma^2 \in \mathcal{G}^2.
\end{align*}
Note that Nash equilibria in a team game are also known as person-by-person optimal solutions.
\end{definition}

\section{Rate-Regulation Tradeoff}\label{c1:problem-formulation}
Consider a Gauss--Markov process with the discrete-time time-varying state equation
\begin{align}\label{c1:eq:sys}
	x_{k+1} &= A_k x_k + B_k u_k + w_k,
\end{align}
for $ k \in \mathbb{N}_{[0,N]}$ with initial condition $x_0$, where $x_k \in \mathbb{R}^n$ is the state of the process, $A_k \in \mathbb{R}^{n \times n}$ is the state matrix, $B_k \in \mathbb{R}^{n \times m}$ is the input matrix, $u_k \in \mathbb{R}^m$ is the control input applied by an actuator and decided by a controller that is collocated with the actuator, $w_k \in \mathbb{R}^n$ is a Gaussian white noise with zero mean and covariance $W_k \succ 0$, and $N \in \mathbb{N}$ is a finite time horizon. A noisy version of the output of the process is observed by a sensor at each time $k$, and given by the output equation
\begin{align}\label{c1:eq:output}
	y_k &= C_k x_k + v_k,
\end{align}
for $k \in \mathbb{N}_{[0,N]}$, where $y_k \in \mathbb{R}^p$ is the output of the process, $C_k \in \mathbb{R}^{p \times n}$ is the output matrix, and $v_k \in \mathbb{R}^p$ is a Gaussian white noise with zero mean and covariance  $V_k \succ 0$. It is assumed that $x_0$ is a Gaussian vector with mean $m_0$ and covariance $M_0$, and that $x_0$, $w_k$, and $v_k$ are mutually independent for all $k \in \mathbb{N}_{[0,N]}$.

The feedback control loop is closed via a reliable but costly communication channel, and the sensory information in this channel is carried in the form of data packets subject to one-step delay. Let $a_k$ and $b_k$ represent the input and the output of the channel at time $k$, respectively. Then, we have
\begin{align}\label{c1:eq:etm1}
b_{k+1} = \left\{
  \begin{array}{l l}
     a_k, & \ \text{if} \ \delta_k =1, \\
     \varnothing, & \ \text{otherwise},
  \end{array} \right.
\end{align}
for $k \in \mathbb{N}_{[0,N]}$ with $b_0 = \varnothing$, where $\delta_k \in \{0,1\}$ is the transmission decision decided by an event trigger that is collocated with the sensor. It is assumed that the data packet that can be transmitted at time $k$ contains the minimum mean-square-error state estimate at the event trigger at time $k$, and that the quantization error is negligible. Clearly, this state estimate condenses all the previous and current outputs of the process, and its transmission is always better than that of the raw output at time $k$.

The event trigger and the controller, as two distributed decision makers, make their decisions at each time $k$ based on their causal information sets, which are given by
\begin{align}
	\mathcal{I}^e_k &:= \Big\{ y_t, b_t, \delta_{s}, u_{s} \Big| t \in \mathbb{N}_{[0,k]}, s \in \mathbb{N}_{[0,k-1]} \Big\},\\[1.45\jot]
	\mathcal{I}^c_k &:= \Big\{ b_t, \delta_{s}, u_{s} \Big| t \in \mathbb{N}_{[0,k]}, s \in \mathbb{N}_{[0,k-1]} \Big\},
\end{align}
respectively, We say that a triggering policy $\pi$ and a control policy $\mu$ are admissible if $\pi = \{\ProbM(\delta_k | \mathcal{I}^e_k) \}_{k=0}^{N}$ and $\mu = \{\ProbM(u_k | \mathcal{I}^c_k) \}_{k=0}^{N}$, where $\ProbM(\delta_k | \mathcal{I}^e_k)$ and $\ProbM(u_k | \mathcal{I}^c_k)$ are Borel measurable stochastic kernels defined on suitable measurable spaces. We represent the sets of admissible triggering policies and admissible control policies by $\mathcal{P}$ and $\mathcal{M}$, respectively. 

For the system outlined above, we are interested in a rate-regulation tradeoff between the packet rate and the regulation cost. Let us measure the packet rate by
\begin{equation}\label{eq:rate-measure}
R(\pi,\mu) := \textstyle \frac{1}{N+1} \E\Big[\sum_{k=0}^{N} \ell_k \delta_k \Big],
\end{equation}
where $\ell_k \geq 0$ is a weighting coefficient, and measure the regulation cost by
\begin{equation}\label{eq:control-measure}
\begin{aligned}
J(\pi,\mu) := \textstyle \frac{1}{N+1} \E \Big [\textstyle\sum_{k=0}^{N+1} x_{k}^T Q_{k} x_{k} + \textstyle\sum_{k=0}^{N} u_k^T R_k u_k \Big],
\end{aligned}
\end{equation}
where $Q_k \succeq 0$ and $R_k \succ 0$ are weighting matrices. The rate-regulation tradeoff can then be expressed as a stochastic optimization problem with the loss function
\begin{align}\label{eq:main_problem1}
	\Phi(\pi,\mu) := (1-\lambda) R(\pi,\mu) + \lambda J(\pi,\mu),
\end{align}
over the space of admissible policy profiles $(\pi,\mu) \in \mathcal{P}\times \mathcal{M}$ given a tradeoff multiplier $\lambda \in (0,1)$. This tradeoff, as we will see, allows us to describe the value of information.

\begin{remark}
The rate-regulation tradeoff, which is formulated based on the weighted sum approach (see e.g., \cite{marler2010}), is a tradeoff between two objective functions. The objective function in (\ref{eq:rate-measure}) penalizes the packet rate in the communication channel, and is appropriate for packet switching networks. This objective function takes into account the price of communication through the weighting coefficient. Moreover, the objective function in (\ref{eq:control-measure}) penalizes the state deviation and the control effort, and is appropriate for regulation tasks. This objective function can be modified for tracking tasks by a transformation when the reference trajectory is known. Finally, note that the underlying optimization problem with the loss function (\ref{eq:main_problem1}) over the space of admissible policy profiles $(\pi,\mu) \in \mathcal{P}\times \mathcal{M}$ is in general an intractable problem. However, in this article, based on a game theoretic analysis, we prove the existence of a Nash equilibrium $(\pi^{\star},\mu^{\star})$. Even though we investigate an imperfect information case, the results can easily be specialized for the perfect information counterpart.
\end{remark}

\section{Quantification of\\the Value of Information}\label{c1:main-results1}
In this section, we present our results on the characterization and the computation of the value of information. We first show how the value of information emerges from the rate-regulation tradeoff formulated in the previous section. We then discuss its structural properties and computational aspects.

\subsection{Formula of the Value of Information}
Since in the rate-regulation tradeoff the system  has two decision makers with different information sets, we can define two different value functions, viz., one from the perspective of the event trigger, i.e., $V^e_k(\mathcal{I}^e_k)$, and one from the perspective of the controller, i.e., $V^c_k(\mathcal{I}^c_k)$. Based on this observation, we introduce our general formula of the value of information in the following definition.

\begin{definition}[Value of Information]\label{def:voi}
The value of information at time $k$ is defined as the variation in the value function $V^e_k(\mathcal{I}^e_k)$ with respect to the sensory information $a_k$ that can be communicated to the controller at time~$k$, i.e.,
\begin{align}\label{eq:voi-def}
\voi_k(\mathcal{I}^e_k) := V^e_k(\mathcal{I}^e_k)|_{\delta_k = 0} - V^e_k(\mathcal{I}^e_k)|_{\delta_k = 1},
\end{align}
where $V^e_k(\mathcal{I}^e_k)|_{\delta_k}$ denotes the value function $V^e_k(\mathcal{I}^e_k)$ when the transmission decision $\delta_k$ is enforced.
\end{definition}

\begin{remark}
The value of information $\voi_k(\mathcal{I}^e_k)$, defined in (\ref{eq:voi-def}), in a sense measures the sensitivity of the value function $V^e_k(\mathcal{I}^e_k)$ with respect to a data packet that can be transmitted to the controller at time $k$. Note that the above formula is general and valid for any choice of the system model. Furthermore, recall that we are interested in a valuation of information associated with a decision about transmission of a data packet through the observation channel at each time $k$. This decision is made by the event trigger and according to the stochastic kernel $\ProbM(\delta_k | \mathcal{I}^e_k)$. For this reason, $\voi_k(\mathcal{I}^e_k)$ was evaluated based on the value function $V^e_k(\mathcal{I}^e_k)$, and not the value function $V^c_k(\mathcal{I}^c_k)$.
\end{remark}

The next lemma introduces a loss function that is equivalent to the original loss function in the sense that it yields the same optimal decision policies. Associated with this loss function, we will subsequently define the value functions $V^e_k(\mathcal{I}^e_k)$ and $V^c_k(\mathcal{I}^c_k)$.

\begin{lemma}\label{lem:1}
Let $S_k \succeq 0$ be a matrix obeying the algebraic Riccati equation
\begin{equation}\label{eq:riccati}
\begin{aligned}
S_k &= Q_k + A_k^T S_{k+1} A_k - A_k^T S_{k+1} B_k\\[2.65\jot]
	&\qquad \quad \times (B_k^T S_{k+1} B_k + R_k)^{-1} B_k^T S_{k+1} A_k,
\end{aligned}
\end{equation}
for $k \in \mathbb{N}_{[0,N]}$ with initial condition $S_{N+1} = Q_{N+1}$. Then,
\begin{align}\label{eq:phiprime}
	\Psi&(\pi,\mu) := \E \Big[ \textstyle \sum_{k=0}^{N} \theta_k \delta_k + \varsigma_k \Big],
\end{align}
is equivalent to $\Phi(\pi,\mu)$, where $\theta_k = \ell_k (1-\lambda)/\lambda$ and $\varsigma_k =  (u_k + (B_k^T S_{k+1} B_k + R_k)^{-1} B_k^T S_{k+1} A_k x_k)^T (B_k^T S_{k+1} B_k + R_k) (u_k + (B_k^T S_{k+1} B_k + R_k)^{-1} B_k^T S_{k+1} A_k x_k)$.
\end{lemma}

\begin{proof}
The result is proved by applying few operations on the state equation (\ref{c1:eq:sys}) and the algebraic Riccati equation (\ref{eq:riccati}), and by discarding the fixed terms (see e.g., \cite{stoccontrol}).
\end{proof}

\begin{definition}[Value functions]
The value functions $V^e_k(\mathcal{I}^e_k)$ and $V^c_k(\mathcal{I}^c_k)$ are defined as
\begin{align}
	V^e_k(\mathcal{I}^e_k) :=& \min_{\pi \in \mathcal{P} : \mu = \mu^\star}\E \Big[ \textstyle \sum_{t=k}^{N} \theta_t \delta_t + \varsigma_{t+1} \Big| \mathcal{I}^e_k \Big],\label{eq:Ve-def}\\[1.5\jot]
	V^c_k(\mathcal{I}^c_k) :=& \min_{\mu \in \mathcal{M}: \pi = \pi^\star}\E \Big[ \textstyle \sum_{t=k}^{N} \theta_{t-1} \delta_{t-1} + \varsigma_{t} \Big| \mathcal{I}^c_k \Big],\label{eq:Vc-def}
\end{align}
for $k \in \mathbb{N}_{[0,N]}$ given a policy profile $(\pi^{\star},\mu^{\star})$, where we adopt the convention $\theta_{-1} =0$, $\varsigma_{N+1} = 0$, and $S_{N+2} = 0$.
\end{definition}

Note that at a Nash equilibrium the value functions $V^e_k(\mathcal{I}^e_k)$ and $V^c_k(\mathcal{I}^c_k)$ should simultaneously satisfy the optimality relations (see e.g., \cite{basargame}).

\subsection{Emergence of the Value of Information}
Let $\check{x}_k$ and $\hat{x}_k$, unless otherwise stated, denote the minimum mean-square-error state estimates at the event trigger and the controller at time $k$, respectively. In addition, let us define the estimation error from the perspective of the event trigger $\check{e}_k := x_k - \check{x}_k$, the estimation error from the perspective of the controller $\hat{e}_k := x_k - \hat{x}_k$, and the estimation mismatch $\tilde{e}_k := \check{x}_k - \hat{x}_k$. The following two propositions characterize the optimal estimators at the event trigger and the controller. These estimators will be needed for our structural result.

\begin{proposition}\label{c1:prop:imperf-KF}
The conditional mean $\E[{x}_k | \mathcal{I}^e_k]$ is the minimum mean-square-error estimator at the event trigger, and obeys
\begin{align}
\begin{split}\label{eq:kf-E}
	\check{x}_{k+1} &= A_k \check{x}_k + B_k u_k\\[2.45\jot]
	&\quad + K_{k+1} \big(y_{k+1} - C_{k+1} ( A_k \check{x}_k + B_k u_k)\big),
\end{split}\\[2.45\jot]
\begin{split}\label{eq:kf-Cov}
	Y_{k+1} &= \big( (A_k Y_k A_k^T + W_k)^{-1} + C_{k+1}^T V_{k+1}^{-1} C_{k+1} \big)^{-1},	
\end{split}
\end{align}
for $k \in \mathbb{N}_{[0,N]}$ with initial conditions $\check{x}_0 = m_0 + Y_{0} C_{0}^T V_{0}^{-1}(y_0 - C_0 m_0)$ and $Y_0 = (M_0^{-1} + C_{0}^T V_{0}^{-1} C_{0})^{-1}$, where $\check{x}_k = \E[{x}_k | \mathcal{I}^e_k]$, $Y_k = \Cov[x_k | \mathcal{I}^e_k]$, and $K_k = Y_k C_k^T V_k^{-1}$.
\end{proposition}

\begin{proof}
Clearly, given the information set $\mathcal{I}^e_k$, the mean-square error is minimized by $\E[{x}_k | \mathcal{I}^e_k]$, and the optimal estimator is the standard Kalman filter (see e.g., \cite{stengel1994}).
\end{proof}

\begin{proposition}\label{c1:prop:imperf-estimatorX}
The conditional mean $\E[{x}_k | \mathcal{I}^c_k]$ is the minimum mean-square-error estimator at the controller, and obeys
\begin{equation}\label{c1:eq:imperf-estimate-dynX}
\begin{aligned}
	\hat{x}_{k+1} &= A_k \hat{x}_k + B_k u_k + \delta_k A_k \tilde{e}_k + (1-\delta_k) \imath_k,
\end{aligned}
\end{equation}
for $k \in \mathbb{N}_{[0,N]}$ with initial condition $\hat{x}_0 = m_0$, where $\hat{x}_k = \E[{x}_k | \mathcal{I}^c_k]$ and $\imath_k = A_k \E[\hat{e}_k | \mathcal{I}^c_k,\delta_k=0]$. In addition, the conditional covariance $\Cov[{x}_k | \mathcal{I}^c_k]$ obeys
\begin{equation}\label{c1:eq:imperf-cov-dynX}
\begin{aligned}
	Z_{k+1} &= A_k Z_k A_k^T + W_k\\[2.05\jot]
	&\qquad \quad \ - \delta_k A_k (Z_k - Y_k) A_k^T - (1-\delta_k) \Xi_k,	
\end{aligned}
\end{equation}
for $k \in \mathbb{N}_{[0,N]}$ with initial condition $Z_0 = M_0$, where $Z_k = \Cov[x_k | \mathcal{I}^c_k]$ and $\Xi_k = A_k (Z_k - \Cov[\hat{e}_k | \mathcal{I}_k^c, \delta_k =0]) A_k^T$.
\end{proposition}

\begin{proof}
Clearly, given the information set $\mathcal{I}^c_k$, the mean-square error is minimized by $\E[{x}_k | \mathcal{I}^c_k]$. Moreover, from the state equation (\ref{c1:eq:sys}), we see that
\begin{align}
\hat{x}_{k+1} &= A_k \E[x_k | \mathcal{I}_{k+1}^c] + B_k u_k,\label{c1:eq:imperf-propag-meanX}\\[2.3\jot]
Z_{k+1} &= A_k \Cov[x_k | \mathcal{I}_{k+1}^c] A_k^T + W_k.\label{c1:eq:imperf-propag-covX}
\end{align}
By definition, the transmission decision $\delta_k$ at each time can be either one or zero. If $\delta_k = 1$, the controller receives $\check{x}_k$ at time $k+1$. In this case, we can write
\begin{align*}
\Prob(x_k | \mathcal{I}^c_{k+1}) = \Prob(x_k | \mathcal{I}_k^c, \check{x}_k, Y_k) = \Prob(x_k | \mathcal{I}^e_k),
\end{align*}
where we used the fact that $\{\check{x}_k,Y_k\}$ is statistically equivalent to $\mathcal{I}^e_k$. Hence, we obtain $\E[x_k | \mathcal{I}_{k+1}^c]  = \check{x}_k$ and $\Cov[x_k | \mathcal{I}_{k+1}^c]  = Y_k$. However, if $\delta_k = 0$, the controller receives nothing at time $k+1$. In this case, we can write
\begin{align*}
\Prob(x_k | \mathcal{I}^c_{k+1}) &= \Prob(x_k | \mathcal{I}_k^c, \delta_k = 0)\\[2.8\jot]
	&= \textstyle \frac{1}{c} \Prob(\delta_k = 0 | \mathcal{I}_k^c, x_k) \Prob(x_k | \mathcal{I}_k^c),
\end{align*}
where $c$ is a normalizer. Hence, for any admissible triggering policy, the conditional mean $\E[x_k | \mathcal{I}_{k+1}^c]$ and the conditional covariance $\Cov[x_k | \mathcal{I}_{k+1}^c]$ can be computed based on $\Prob(x_k | \mathcal{I}^c_{k+1})$. Let us define $x'_k$ and $Z'_k$ as
\begin{align*}
x'_k &:= \E[x_k | \mathcal{I}_k^c, \delta_k =0] - \hat{x}_k,\\[2.5\jot]
Z'_k &:= Z_k - \Cov[x_k | \mathcal{I}_k^c, \delta_k =0].
\end{align*}
Consequently, for any value of $\delta_k$, we see that
\begin{align}
\E[x_k | \mathcal{I}_{k+1}^c]  &= \hat{x}_k + \delta_k (\check{x}_k - \hat{x}_k) + (1-\delta_k) x'_k,\label{c1:eq:imperf-update-meanX}\\[2.4\jot]
\Cov[x_k | \mathcal{I}_{k+1}^c]  &= Z_k - \delta_k (Z_k - Y_k) - (1-\delta_k) Z'_k.\label{c1:eq:imperf-update-covX}
\end{align}
Now, we only need to substitute~(\ref{c1:eq:imperf-update-meanX}) and (\ref{c1:eq:imperf-update-covX}) in~(\ref{c1:eq:imperf-propag-meanX}) and (\ref{c1:eq:imperf-propag-covX}), respectively, and define $\imath_k := A_k x'_k$ and $\Xi_k := A_k Z'_k A_k^T$. We can write $\imath_k = A_k ( \E[ x_k | \mathcal{I}^c_k, \delta_k =0 ] - \E[ x_k | \mathcal{I}^c_k]) = A_k \E[ \hat{e}_k | \mathcal{I}_k^c, \delta_k =0]$ and $\Xi_k = A_k (Z_k - \Cov[x_k | \mathcal{I}^c_k, \delta_k=0]) A_k^T = A_k (Z_k - \Cov[\hat{e}_k | \mathcal{I}^c_k, \delta_k=0]) A_k^T$, and the result follows.
\end{proof}

\begin{remark}
Observe that the optimal estimators at the event trigger and the controller have completely different structures. While the conditional distribution $\ProbM(x_k | \mathcal{I}^e_k)$ is Gaussian and the conditional mean $\check{x}_k$ obeys a linear recursive equation, the conditional distribution $\ProbM(x_k | \mathcal{I}^c_k)$ is generally non-Gaussian and the conditional mean $\hat{x}_k$ generally obeys a nonlinear recursive equation. Note that the residuals $\imath_k$ and $\Xi_k$ in (\ref{c1:eq:imperf-estimate-dynX}) and (\ref{c1:eq:imperf-cov-dynX}) are both due to negative information. The existence of these terms implies that the controller might be able to decrease its uncertainty even when it does not receive any data packet from the sensor. The values of the residuals $\imath_k$ and $\Xi_k$ at each time $k$ depend on the structure of the triggering policy. For any fixed triggering policy, these values can be computed numerically by techniques from nonlinear filtering (see e.g., \cite{sijs2012,wu2013}).
\end{remark}

The next theorem, which is our main result, characterizes a Nash equilibrium in the rate-regulation tradeoff.

\begin{theorem}\label{c1:thm:imperf-optimalityX}
There exists at least one Nash equilibrium $(\pi^{\star},\mu^{\star})$ in the rate-regulation tradeoff such that 
\begin{align}\label{eq:opt-profile}
(\pi^{\star},\mu^{\star}) = \Big( \big\{\mathds{1}_{\voi_k \geq 0} \big\}_{k=0}^{N}, \big\{ - L_k \hat{x}_k \big\}_{k=0}^{N} \Big),
\end{align}
with the value of information $\voi_k(\mathcal{I}^e_k)$ as a symmetric function of the estimation mismatch $\tilde{e}_k$ obeying
\begin{align}
	\voi_k(\mathcal{I}^e_k) &= \tilde{e}_k^T A_k^T \Gamma_{k+1} A_k \tilde{e}_k - \theta_k + \varrho_k,\label{eq:voi-imperfectX}
\end{align}
where $\Gamma_k = A_k^T S_{k+1} B_k (B_k^T S_{k+1} B_k + R_k)^{-1} B_k^T S_{k+1} A_k$ and $\varrho_k = \E[V^e_{k+1}(\mathcal{I}^e_{k+1})|\mathcal{I}^e_k, \delta_k = 0] - \E[V^e_{k+1}(\mathcal{I}^e_{k+1})|\mathcal{I}^e_k, \delta_k = 1]$, and with the conditional mean $\hat{x}_k$ without being affected by negative information obeying
\begin{align}
\hat{x}_{k+1} &= A_k \hat{x}_k + B_k u_k + \delta_k A_k \tilde{e}_k,
\end{align}
for $k \in \mathbb{N}_{[0,N]}$ with initial condition $\hat{x}_0 = m_0$, where $L_k = (B_k^T S_{k+1} B_k + R_k)^{-1} B_k^T S_{k+1} A_k$.
\end{theorem}

\begin{proof}
The proof is structured in two parts. In the first part, we show that $\Psi(\pi^{\star},\mu^{\star}) \leq \Psi(\pi,\mu^{\star})$ for all $\pi \in \mathcal{P}$. Note that given the control policy $\mu^\star$, the state estimate $\hat{x}_k$ obeys $\hat{x}_{k+1} = A_k \hat{x}_k + B_k u_k + \delta_k A_k \tilde{e}_k$ for $k \in \mathbb{N}_{[0,N]}$ with initial condition $\hat{x}_0 = m_0$. From the additivity of $V^e_k(\mathcal{I}^e_k)$, we obtain
\begin{align*}
	V^e_k(\mathcal{I}^e_k) &= \min_{\Prob(\delta_k|\mathcal{I}_k^e)} \E\Big[ \theta_k \delta_k + \hat{e}_{k+1}^T \Gamma_{k+1} \hat{e}_{k+1} \\[-0.25\jot]
	&\quad +\min_{\Prob(\delta_{k+1}|\mathcal{I}_{k+1}^e)} \E\Big[ \theta_{k+1} \delta_{k+1}\\[0.25\jot]
	&\quad +  \hat{e}_{k+2}^T \Gamma_{k+2} \hat{e}_{k+2} + \dots \Big|\mathcal{I}^e_{k+1} \Big] \Big| \mathcal{I}^e_k\Big]\\[1.75\jot]
	&=\min_{\Prob(\delta_k|\mathcal{I}_k^e)} \E\Big[\theta_k \delta_k + \hat{e}_{k+1}^T \Gamma_{k+1} \hat{e}_{k+1} + V^e_{k+1}(\mathcal{I}^e_{k+1}) \Big|\mathcal{I}^e_k\Big],
\end{align*}
for $k \in \mathbb{N}_{[0,N]}$ with initial condition $V^e_{N+1}(\mathcal{I}^e_{N+1}) = 0$. We prove by backward induction that $V^e_k(\mathcal{I}^e_k)$ is a symmetric function of $\tilde{e}_k$. Clearly, the claim is satisfied for time $N+1$. We assume that the claim holds at time $k+1$, and shall prove that it also holds at time $k$. Given the dynamics of $\hat{x}_k$ in this case, we observe that $\hat{e}_k$ and $\tilde{e}_k$ obey
\begin{align}
	\hat{e}_{k+1} &= A_k \hat{e}_k - \delta_k A_k \tilde{e}_k + w_k,\label{eq:prove-error-dyn-imperfX}\\[2.8\jot]
	\tilde{e}_{k+1} &= (1-\delta_k) A_k \tilde{e}_k + n_k,\label{eq:xi-dynX}
\end{align}
for $k \in \mathbb{N}_{[0,N]}$ with initial conditions $\hat{e}_0 = x_0 - \hat{x}_0$ and $\tilde{e}_0 = \check{x}_0 - \hat{x}_0$, where $n_k \in \mathbb{R}^n$ is a Gaussian white noise with zero mean and covariance $N_k = K_{k+1} ( C_{k+1}(A_k Y_k A_k^T + W_k) C_{k+1}^T + V_{k+1}) K_{k+1}^T$. From (\ref{eq:prove-error-dyn-imperfX}), we find
\begin{align*}
	&\E\Big[ \hat{e}_{k+1}^T \Gamma_{k+1} \hat{e}_{k+1} \Big| \mathcal{I}^e_k \Big]\\[1.85\jot]
	&\quad = \E \Big[ \hat{e}_k^T A_k^T \Gamma_{k+1} A_k \hat{e}_k + \delta_k^2 \tilde{e}_k^T A_k^T \Gamma_{k+1} A_k \tilde{e}_k\\[2.35\jot]
	&\quad \ \qquad + w_k^T \Gamma_{k+1} w_k - 2 \delta_k \tilde{e}_k^T A_k^T \Gamma_{k+1} A_k \hat{e}_k\\[2.35\jot]
	&\quad \ \qquad - 2 \delta_k \tilde{e}_k^T A_k^T \Gamma_{k+1} w_k + 2 \hat{e}_k^T A_k^T \Gamma_{k+1} w_k \Big| \mathcal{I}^e_k \Big]\\[1.8\jot]
	&\quad\ = \underset{\delta_k}{\E} \Big[ \tilde{e}_k^T A_k^T \Gamma_{k+1} A_k \tilde{e}_k + \tr(A_k^T \Gamma_{k+1} A_k Y_k)\\[.85\jot]
	&\quad \ \qquad +\tr(\Gamma_{k+1} W_k) - \delta_k \tilde{e}_k^T A_k^T \Gamma_{k+1} A_k \tilde{e}_k \Big| \mathcal{I}_k^e \Big],
\end{align*}
where we used the facts that $\E[\hat{e}_k | \mathcal{I}^e_k] = \tilde{e}_k$, $\Cov[\hat{e}_k | \mathcal{I}^e_k]=Y_k$, $\E[\tilde{e}_k | \mathcal{I}^e_k] = \tilde{e}_k$, $\E[w_k | \mathcal{I}^e_k] =0$, and that $w_k$ is independent of $\hat{e}_k$. Accordingly, we can show that
\begin{equation}\label{c1:eq:cost-to-go2-imperfectX}
\begin{aligned}
	&V^e_k(\mathcal{I}^e_k) = \min_{\delta_k} \Big\{\theta_k \delta_k + (1-\delta_k) \tilde{e}_k^T A_k^T \Gamma_{k+1} A_k \tilde{e}_k\\[1.25\jot]
	&\ + \tr(\Gamma_{k+1} W_k + A_k^T \Gamma_{k+1} A_k Y_k) + \E[V^e_{k+1}(\mathcal{I}^e_{k+1})|\mathcal{I}^e_k]  \Big\}.
\end{aligned}
\end{equation}
The minimizer in (\ref{c1:eq:cost-to-go2-imperfectX}) is obtained as \vspace{-1mm} $\delta_k^{\star} = \mathds{1}_{\voi_k(\mathcal{I}^e_k) \geq 0}$, where 
\begin{align*}
\voi_k(\mathcal{I}^e_k) = \tilde{e}_k^T A_k^T \Gamma_{k+1} A_k \tilde{e}_k - \theta_k + \varrho_k,	
\end{align*}
and $\varrho_k = \E[V^e_{k+1}(\mathcal{I}^e_{k+1})|\mathcal{I}^e_k, \delta_k = 0] - \E[V^e_{k+1}(\mathcal{I}^e_{k+1})|\mathcal{I}^e_k, \delta_k = 1]$. Define now $\bar{n}_k := - n_k$. Note that $\bar{n}_k$ is also a Gaussian white noise with zero mean and covariance $N_k$. It follows that
\begin{align*}
	\E &\Big[V^e_{k+1}( \tilde{e}_{k+1}) \Big|\mathcal{I}^e_k, \delta_k \Big]\\[1.85\jot]
	&= \E \Big[V^e_{k+1} \big( (1-\delta_k) A_k \tilde{e}_k + n_k \big) \Big|\mathcal{I}^e_k, \delta_k \Big]\\[2\jot]
	&= \E \Big[V^e_{k+1} \big(- (1-\delta_k) A_k \tilde{e}_k - n_k \big) \Big|\mathcal{I}^e_k, \delta_k \Big]\\[2\jot]
	&= \textstyle \int_{\mathbb{R}^n} V^e_{k+1} \big(-(1-\delta_k) A_k \tilde{e}_k - n_k \big) d \exp(-\frac{1}{2} n_k^T N_k^{-1} n_k) d n_k\\[3\jot]
	&= \textstyle \int_{\mathbb{R}^n} V^e_{k+1} \big(-(1-\delta_k) A \tilde{e}_k + \bar{n}_k \big) d  \exp(-\frac{1}{2} \bar{n}_k^T N_k^{-1} \bar{n}_k) d \bar{n}_k\\[2.7\jot]	
	&= \E \Big[V^e_{k+1} \big(- (1-\delta_k) A_k \tilde{e}_k + n_k \big) \Big|\mathcal{I}^e_k, \delta_k \Big],
\end{align*}
where $d$ is a constant, the first equality comes from (\ref{eq:xi-dynX}), and the second equality from the hypothesis assumption. Therefore, $\E[V^e_{k+1}(\mathcal{I}^e_{k+1}) |\mathcal{I}^e_k, \delta_k]$ is a symmetric function of $\tilde{e}_k$. This implies that $\varrho_k$ and $\voi_k(\mathcal{I}^e_k)$ are also symmetric functions of $\tilde{e}_k$. Moreover, we can write $V^e_k(\mathcal{I}^e_k)$~as
\begin{align}
V^e_k(\mathcal{I}^e_k) = \left\{
\begin{array}{l l}
     V^{e1}_k(\mathcal{I}^e_k), & \ \text{if} \ \voi_k(\mathcal{I}^e_k) \geq 0, \\[1.5\jot]
     V^{e0}_k(\mathcal{I}^e_k), & \ \text{otherwise},
\end{array} \right.
\end{align}
where $V^{e1}_k(\mathcal{I}^e_k)$ and $V^{e0}_k(\mathcal{I}^e_k)$ are both symmetric functions of $\tilde{e}_k$. Hence, we conclude that $V^e_k(\mathcal{I}^e_k)$ is a symmetric function of $\tilde{e}_k$.

In the second part, we show that $\Psi(\pi^{\star},\mu^{\star}) \leq \Psi(\pi^{\star},\mu)$ for all $\mu \in \mathcal{M}$. Note that given the triggering policy $\pi^\star$, the state estimate $\hat{x}_k$ obeys $\hat{x}_{k+1} = A_k \hat{x}_k + B_k u_k + \delta_k A_k \tilde{e}_k + (1-\delta_k) \imath_k$ for all $k \in \mathbb{N}_{[0,N]}$ with initial condition $\hat{x}_0 = m_0$. From the additivity of $V^c_k(\mathcal{I}^c_k)$, we obtain
\begin{align*}
	V^c_k(\mathcal{I}^c_k) &= \min_{\Prob(u_k|\mathcal{I}_k^c)} \E \Big[ c_{k-1} + (u_k + L_k x_k)^T \Lambda_k (u_k + L_k x_k) \\[-.5\jot]
	&\quad + \min_{\Prob(u_{k+1}|\mathcal{I}_{k+1}^c)} \E \Big[  c_k + (u_{k+1} + L_{k+1} x_{k+1})^T \Lambda_{k+1}\\[0\jot]
	&\quad \times (u_{k+1} + L_{k+1} x_{k+1}) + \dots \Big| \mathcal{I}^c_{k+1}\Big] \Big| \mathcal{I}^c_k \Big]\\[1.2\jot]
	&= \min_{\Prob(u_k|\mathcal{I}_k^c)} \E \Big[ c_{k-1} + (u_k + L_k x_k)^T \Lambda_k\\[-.5\jot]
	&\qquad \qquad \qquad \qquad \times (u_k + L_k x_k) + V^c_{k+1}(\mathcal{I}^c_{k+1}) \Big| \mathcal{I}^c_k\Big],
\end{align*}
for $k \in \mathbb{N}_{[0,N]}$ with initial condition $V^c_{N+1}(\mathcal{I}^c_{N+1}) = 0$, where $c_k = \theta_k \mathds{1}_{\voi_k(\mathcal{I}^e_k) \geq 0}$ is a function of $\tilde{e}_k$ and $\Lambda_k = B_k^T S_{k+1} B_k + R_k$. We prove by backward induction that $V^c_k(\mathcal{I}^c_k)$ is independent of the previous control inputs. Clearly, the claim is satisfied for time $N+1$ We assume that the claim holds at time $k+1$, and shall prove that it also holds at time $k$. Using the identity $x_k = \hat{x}_k + \hat{e}_k$, we find 
\begin{align*}
	&\E \Big[ (u_k + L_k x_k)^T \Lambda_k (u_k + L_k x_k) \Big| \mathcal{I}^c_k \Big]\\[1.65\jot]
	&\ = \E\Big[(u_k + L_k \hat{x}_k)^T \Lambda_k (u_k + L_k \hat{x}_k) + \hat{e}_k^T L_k^T \Lambda_k L_k \hat{e}_k\\[1.55\jot]
	&\ \qquad + 2 (u_k + L_k \hat{x}_k)^T \Lambda_k L_k \hat{e}_k  \Big| \mathcal{I}^c_k \Big]\\[1.5\jot]
	&\ = \underset{u_k}{\E} \Big[ (u_k + L_k \hat{x}_k)^T \Lambda_k (u_k + L_k \hat{x}_k) + \tr(\Gamma_k Z_k) \Big| \mathcal{I}_k^c \Big],
\end{align*}
where we used the facts that $\E[\hat{x}_k | \mathcal{I}^c_k] = \hat{x}_k$ and $\E[ \hat{e}_k | \mathcal{I}^c_k] = 0$. Given the dynamics of $\hat{x}_k$ in this case, we observe that $\hat{e}_k$ and $\tilde{e}_k$ obey
\begin{align}
	\hat{e}_{k+1} &= A_k \hat{e}_k - \delta_k A_k \tilde{e}_k + w_k - (1-\delta_k) \imath_k,\\[2.8\jot]
	\tilde{e}_{k+1} &= (1-\delta_k) A_k \tilde{e}_k + n_k - (1-\delta_k) \imath_k,
\end{align}
for $k \in \mathbb{N}_{[0,N]}$ with initial conditions $\hat{e}_0 = x_0 - \hat{x}_0$ and $\tilde{e}_0 = \check{x}_0 - \hat{x}_0$. Since $\delta_k$ is a function of $\tilde{e}_k$, we recursively infer that $\hat{e}_k$ and $\tilde{e}_k$ are independent of the control inputs. Accordingly, we can show that
\begin{equation}
\begin{aligned}\label{c1:eq:cost-to-go1-perfect}
	&V^c_k(\mathcal{I}^c_k) = \min_{u_k} \Big\{ \E[ c_{k-1} | \mathcal{I}^c_k ] + (u_k + L_k \hat{x}_k)^T \Lambda_k \\[1.25\jot]
	&\quad \times (u_k + L_k \hat{x}_k) + \tr(\Gamma_k Z_k) + \E[V^c_{k+1}(\mathcal{I}^c_{k+1}) | \mathcal{I}^c_k] \Big\},
\end{aligned}
\end{equation}
where $\E[ c_{k-1} | \mathcal{I}^c_k]$ and $Z_k = \Cov[ \hat{e}_k | \mathcal{I}^c_k]$ are independent of the control inputs because $\tilde{e}_{k-1}$ and $\hat{e}_k$ are independent of the control inputs, respectively. The minimizer in (\ref{c1:eq:cost-to-go1-perfect}) is obtained as $u^{\star}_k = -L_k \hat{x}_k$, and we conclude that $V^c_k(\mathcal{I}^c_k)$ is independent of the previous control inputs. We now need to prove that $\imath_k = 0$ for all $k \in \mathbb{N}_{[0,N]}$. Note that $\hat{e}_0$ and $\tilde{e}_0$ are Gaussian vectors with zero mean. We assume that $\imath_t = 0$ for all $t \in \mathbb{N}_{[0,k-1]}$, and shall show that $\imath_k = 0$. For any value of $\imath_k$, we have
\begin{align*}
	\Prob(\tilde{e}_k | \mathcal{I}_k^c, \delta_k = 0) \propto \Prob(\delta_k = 0 |\tilde{e}_k)\Prob(\tilde{e}_k | \mathcal{I}_k^c).
\end{align*}
By the hypothesis assumption and using the triggering policy $\pi^{\star}$, we see that
$\Prob(\tilde{e}_k | \mathcal{I}_k^c)$ and $\Prob(\delta_k = 0 | \tilde{e}_k)$  are symmetric with respect to $\tilde{e}_k$. Hence, $\Prob(\tilde{e}_k | \mathcal{I}_k^c, \delta_k = 0)$ is also symmetric with respect to $\tilde{e}_k$. This implies that $\imath_k = 0$, and the proof is complete.
\end{proof}

\begin{remark}
Our structural result shows that at the equilibrium $(\pi^\star,\mu^\star)$ the design of the event trigger and the controller in (\ref{eq:opt-profile}) becomes separated, the optimal estimator at the controller in (\ref{c1:eq:imperf-estimate-dynX}) becomes linear, and the conditional covariance in (\ref{c1:eq:imperf-cov-dynX}) becomes independent of the previous control inputs, implying that the control has no dual effect. In addition, our result shows that $\voi_k(\mathcal{I}^e_k)$ is a symmetric function of the estimation mismatch $\tilde{e}_k$, and that it can be computed with arbitrary accuracy through solving the optimality equation in (\ref{c1:eq:cost-to-go2-imperfectX}) recursively and backward in time. The complexity of this computation is $\mathcal{O}(N d^n s)$ when the estimation mismatch $\tilde{e}_k$ is discretized in a grid with $d^{n}$ points and the associated expected value is obtained based on a weighted sum of $s$ samples.
\end{remark}

\begin{remark}
We argue that instead of fixing an ad-hoc triggering condition and studying the properties of the resulting event-triggered system, i.e., the procedure that has been used in most of the studies on event-triggered estimation and control, one should study a cost-benefit analysis without any limiting assumptions on the information structure or the policy structure, and find a triggering condition as a result of this analysis. Note that $\voi_k(\mathcal{I}^e_k)$ is in fact the difference between the benefit and the cost of a data packet. In light of our structural result, at each time $k$, the benefit of transmitting a data packet is $\tilde{e}_k^T A_k^T \Gamma_{k+1} A_k \tilde{e}_k + \varrho_k$ and its associated cost is $\theta_k$. In this respect, our triggering condition has an important interpretation: a data packet containing the sensory information $\check{x}_k$ should be transmitted to the controller only if its benefit surpasses its cost, i.e., $\voi_k(\mathcal{I}^e_k) \geq 0$. This interpretation does not exist for any triggering condition that is not based on a cost-benefit analysis.
\end{remark}

\begin{remark}
Note that the rate-regulation tradeoff in our study might admit multiple Nash equilibria. Unfortunately, there exists no general procedure for finding all these equilibria (if any). Using  backward induction, we here proved the existence of a Nash equilibrium $(\pi^\star,\mu^\star)$, which has desirable characteristics. Our result guarantees that the set of globally optimal solutions cannot be empty. A natural question that arises in relation to the equilibrium $(\pi^\star,\mu^\star)$ is whether it is globally optimal. We can infer from the results in the literature (see e.g., \cite{lipsa2011,molin2017}) that for the special case of scalar Gauss--Markov processes the optimality gap of this equilibrium is zero. We study this issue for the general case of multi-dimensional Gauss--Markov processes in \cite{voi2}, where we show that the optimality gap of this equilibrium in fact remains zero (see Theorem~1 in \cite{voi2}).
\end{remark}

\subsection{Quadratic Approximation of the Value of Information}
The computation of $\voi_k(\mathcal{I}^e_k)$ based on the optimality equation (\ref{c1:eq:cost-to-go2-imperfectX}) can be difficult especially when $n$ increases. This motivates us to search for an approximation of the value of information that can be expressed analytically. The next proposition provides such an approximation with a performance guarantee.

\begin{proposition}\label{c1:prop-Delta-imperfX}
Let the control policy $\mu^{\star}$ be fixed. A triggering policy $\pi^+$ that outperforms the periodic triggering policy with period one in the rate-regulation tradeoff is given by
\begin{align}
\delta_k^+ = \mathds{1}_{\voi^{+}_k(\mathcal{I}^e_k) \geq 0},
\end{align}
where $\voi^{+}_k(\mathcal{I}^e_k)$ is a quadratic approximation of the value of information expressed as
\begin{align}\label{eq:voi-approx-2}
	\voi^{+}_k(\mathcal{I}^e_k) = \tilde{e}_k^T A_k^T \Gamma_{k+1} A_k \tilde{e}_k - \theta_k.
\end{align}
\begin{proof}
Let $\bar{\pi}$ denote the periodic triggering policy with period one, and $\pi^+$ denote a triggering policy obtained according to
\begin{equation}\label{c1:approx2}
\begin{aligned}
	\delta_k^+ &= \underset{\delta_k}{\argmin} \E \Big[ \theta_k \delta_k + \hat{e}_{k+1}^T \Gamma_{k+1} \hat{e}_{k+1} + V^{\bar{\pi}}_{k+1}(\mathcal{I}^e_{k+1}) \Big|\mathcal{I}^e_k\Big],
\end{aligned}
\end{equation}
where $V^{\bar{\pi}}_k(\mathcal{I}^e_k)$ is the cost-to-go associated with the policy profile $(\bar{\pi},\mu^{\star})$. We prove that $\Psi(\pi^+,\mu^{\star}) \leq \Psi(\bar{\pi},\mu^{\star})$. To do so, it is enough to show $V^{\pi^+}_k \hspace{-1mm}(\mathcal{I}^e_k) \leq V^{\bar{\pi}}_k(\mathcal{I}^e_k)$, where $V^{\pi^+}_k \hspace{-1mm}(\mathcal{I}^e_k)$ is the cost-to-go associated with the policy profile $(\pi^+,\mu^{\star})$. Clearly, $V^{\pi^+}_{N+1}(\mathcal{I}^e_{N+1}) = V^{\bar{\pi}}_{N+1}(\mathcal{I}^e_{N+1}) = 0$. Assume that the claim holds at time $k+1$. We have
\begin{align*}
V^{\pi^+}_k \hspace{-1mm}(\mathcal{I}^e_k) &= \E \Big[ \theta_k \delta_k^{+} + \hat{e}_{k+1}^T \Gamma_{k+1} \hat{e}_{k+1} + V^{\pi^+}_{k+1}(\mathcal{I}^e_{k+1}) \Big|\mathcal{I}^e_k\Big]\\[1.5\jot]
	&\leq \E \Big[ \theta_k \delta_k^{+} + \hat{e}_{k+1}^T \Gamma_{k+1} \hat{e}_{k+1} + V^{\bar{\pi}}_{k+1}(\mathcal{I}^e_{k+1}) \Big|\mathcal{I}^e_k\Big]\\[1.5\jot]
	&\leq \E \Big[ \theta_k + \hat{e}_{k+1}^T \Gamma_{k+1} \hat{e}_{k+1}+ V^{\bar{\pi}}_{k+1}(\mathcal{I}^e_{k+1}) \Big|\mathcal{I}^e_k\Big]\\[2\jot]
	&= V^{\bar{\pi}}_k(\mathcal{I}^e_k),
\end{align*}
where the first inequality comes from the induction hypothesis and the second inequality from the definition of the triggering policy $\pi^+$. Therefore, the claim holds at time $k$.

Besides, following our analysis in the proof of Theorem~\ref{c1:thm:imperf-optimalityX}, we deduce that the minimizer in (\ref{c1:approx2}) is obtained as $\delta_k^{+} = \mathds{1}_{\voi^{+}_k(\mathcal{I}^e_k) \geq 0}$, where
\begin{align*}
	\voi^{+}_k(\mathcal{I}^e_k) = \tilde{e}_k^T A_k^T \Gamma_{k+1} A_k \tilde{e}_k - \theta_k + \varrho^{\bar{\pi}}_k,
\end{align*}
and $\varrho^{\bar{\pi}}_k = \E[V^{\bar{\pi}}_{k+1}(\mathcal{I}^e_{k+1}) | \mathcal{I}_k^e, \delta_k = 0] - \E[V^{\bar{\pi}}_{k+1}(\mathcal{I}^e_{k+1}) | \mathcal{I}_k^e, \delta_k = 1]$. Finally, observe that using the periodic triggering policy $\bar{\pi}$, the estimation error $\hat{e}_{t}$ obeys $\hat{e}_{t+1} = A_t \check{e}_t + w_t$ for $t \in \mathbb{N}_{[k+1,N]}$. Hence, $\hat{e}_{t}$ for $t \in \mathbb{N}_{[k+2,N]}$ is independent of $\delta_k$. This implies that $\varrho^{\bar{\pi}}_k =  0$, and the proof is complete.
\end{proof}
\end{proposition}

\begin{remark}
The value of information approximate $\voi^{+}_k(\mathcal{I}^e_k)$ is a closed-form quadratic function of the estimation mismatch $\tilde{e}_k$, which does not depend on the cost-to-go terms. Our result provides a performance guarantee for this approximation in the sense that the triggering policy $\pi^+$ synthesized based on $\voi^{+}_k(\mathcal{I}^e_k)$ outperforms the periodic triggering policy with period one, when the certainty-equivalent control policy $\mu^\star$ is used. The result is obtained by exploiting a rollout algorithm, which can be viewed as a single iteration of the method of policy iteration.
\end{remark}

\section{Numerical Examples}\label{c1:examples}
\begin{figure}[t]	
\center
\vspace{-2.5mm}
  \includegraphics[width= 0.83\linewidth]{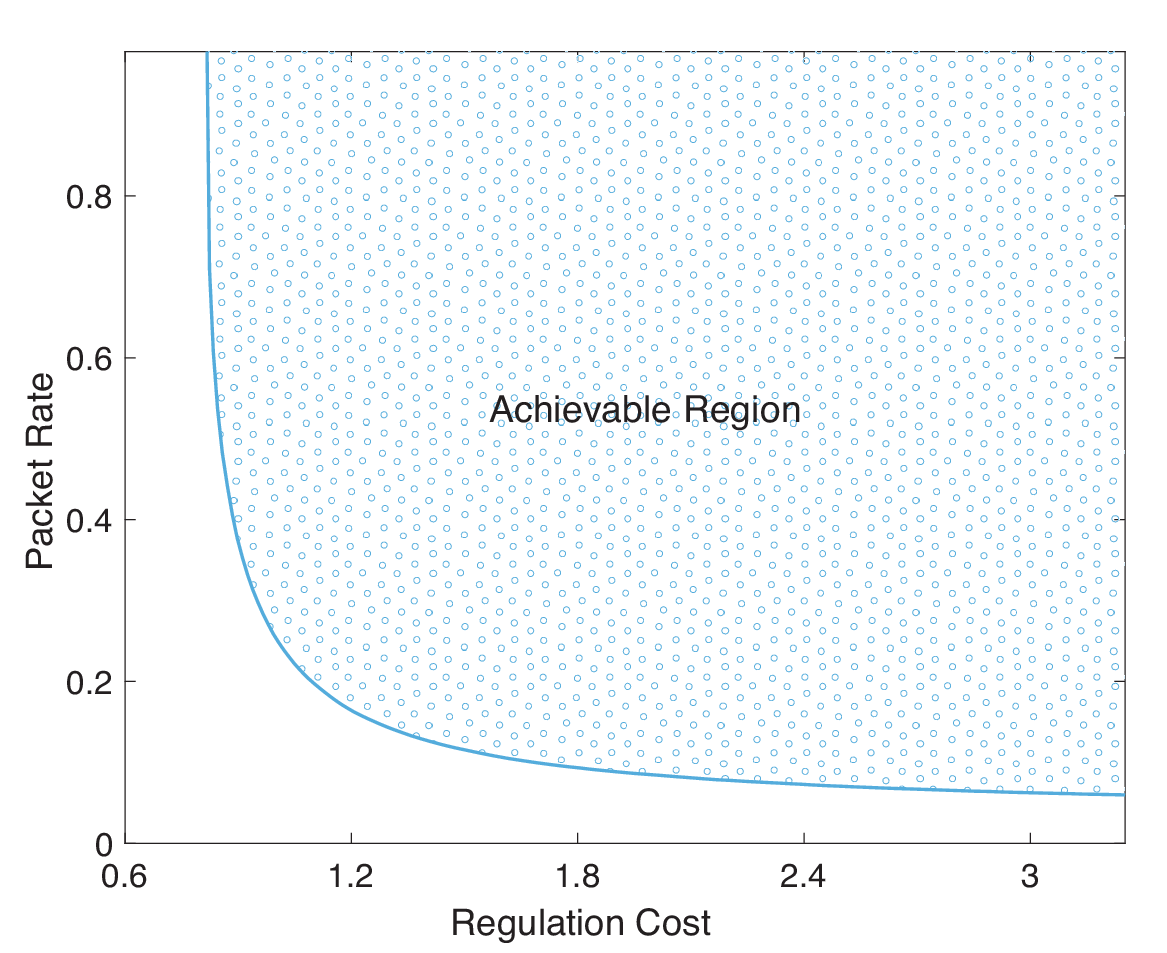}
  \caption{The rate-regulation tradeoff curve. The regulation cost is scaled by one tenth. The achievable region is specified as the area above the tradeoff curve.}
  \label{c1:fig:tradeoff2}
\end{figure}

In this section, we provide two numerical examples that can demonstrate our theoretical results. In the first example, we consider a simple system with state coefficient $A_k = 1.1$, input coefficient $B_k = 1$, output coefficient $C_k = 1$, noise variances $W_k = 3$ and $V_k = 1$ for $k \in \mathbb{N}_{[0,N]}$, mean and variance of the initial condition $m_0 = 0$ and $M_0 = 1$, weighting coefficients $Q_{N+1} = 1$, $\ell_k = 1$, $Q_k = 1$, and $R_k = 0.1$ for $k \in \mathbb{N}_{[0,N]}$, and time horizon $N = 100$. For this system, the rate-regulation tradeoff curve was numerically computed using different values of the tradeoff multiplier $\lambda \in (0,1)$, and is depicted in Fig.~\ref{c1:fig:tradeoff2}. In light of the results in \cite{voi2}, this tradeoff curve is in fact globally optimal. The achievable region is specified in Fig.~\ref{c1:fig:tradeoff2} as the area above the tradeoff curve. Note that there exists no policy profile with performance outside the achievable region.

\begin{figure}[t]
\center
\vspace{-2.5mm}
  \includegraphics[width= 0.83\linewidth]{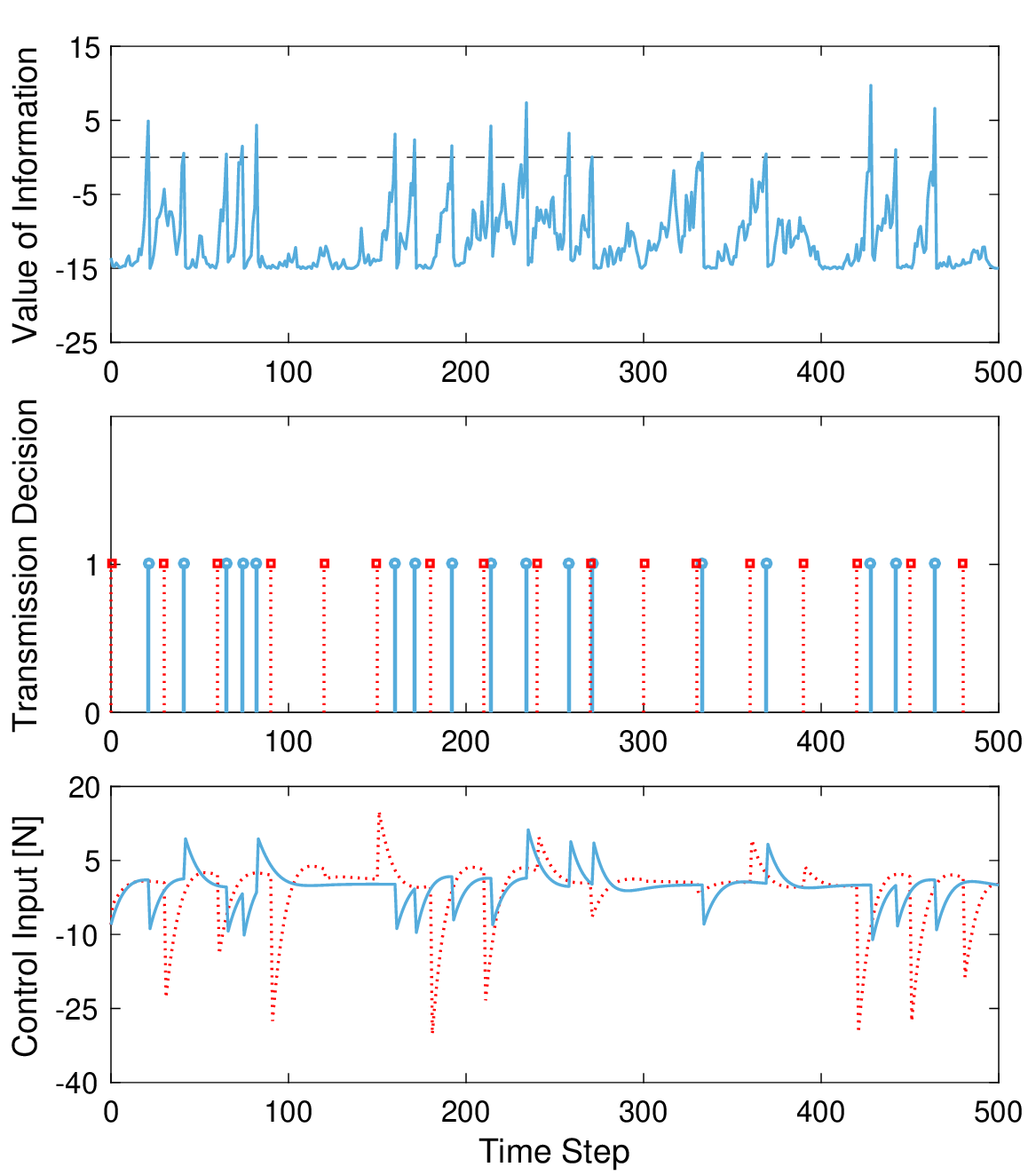}
  \caption{The value of information, transmission decision, and control input trajectories. The value of information is scaled by one tenth. The solid lines represent the trajectories under the triggering policy designed based on the value of information, and the dotted lines represent the trajectories under a periodic triggering policy.}
  \label{c1:fig:trajectories2X}
\end{figure}

\begin{figure}[t]	
\center
\vspace{-2.5mm}
  \includegraphics[width= 0.83\linewidth]{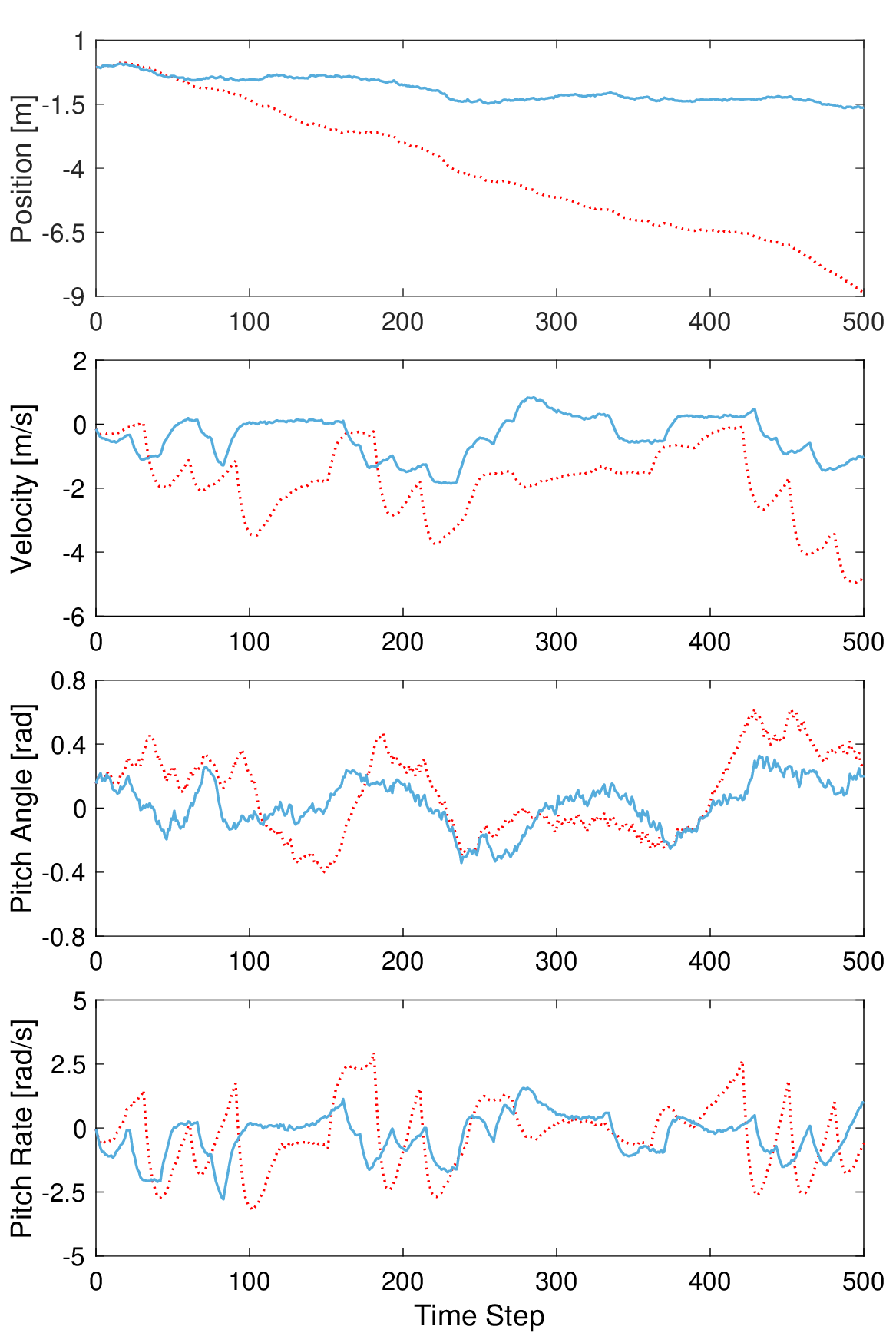}
  \caption{The position, velocity, pitch angle, and pitch rate trajectories. The solid lines represent the trajectories under the triggering policy designed based on the value of information, and the dotted lines represent the trajectories under a periodic triggering policy.}
  \label{c1:fig:trajectories1X}
\end{figure}

In the second example, we consider an inverted pendulum on a cart, for which the continuous-time equations of motion linearized around the unstable equilibrium are given~by
\begin{align*}
	(M+m) \ddot{x} + b \dot{x} - m l \ddot{\phi} = u,\\[2.25\jot]
	(I + m l^2) \ddot{\phi} - m g l \phi = m l \ddot{x},
\end{align*}
where $x$ is the position of the cart, $\phi$ is the pitch angle of the pendulum, $u$ is the force applied to the cart, $M = 0.5 \ \text{kg}$ is the mass of the cart, $m = 0.2 \ \text{kg}$ is the mass of the pendulum, $b = 0.1 \ \text{N/m/sec}$ is the coefficient of friction for the cart, $l = 0.3 \ \text{m}$ is the distance from the pivot to the pendulum's center of mass, $I = 0.006 \ \text{kg.m$^2$}$ is the moment of inertia of the pendulum, and $g = 9.81 \ \text{m/s$^2$}$ is the gravity. We suppose that a sensor measures the position and the pitch angle at each time. The discrete-time state equation of the form~(\ref{c1:eq:sys}), the output equation of the form (\ref{c1:eq:output}), and the loss function of the form (\ref{eq:main_problem1}) are specified~with state, input, and output matrices and noise covariances
\begin{align*}
A_k &= \setlength\arraycolsep{3pt} \begin{bmatrix}
       1.0000    & 0.0100    & 0.0001   & 0.0000\\
       0.0000    & 0.9982    & 0.0267   & 0.0001\\
       0.0000    & 0.0000    & 1.0016   & 0.0100\\
       0.0000    &-0.0045    & 0.3122   & 1.0016
     \end{bmatrix}\!,
B_k = \begin{bmatrix}
        0.0001\\
    	0.0182\\
    	0.0002\\
    	0.0454
     \end{bmatrix}\!,\\[1.5\jot] 
W_k &= \setlength\arraycolsep{3pt} \begin{bmatrix}
    0.0006    & 0.0003    & 0.0001    & 0.0006\\
    0.0003    & 0.0008    & 0.0003    & 0.0004\\
    0.0001    & 0.0003    & 0.0007    & 0.0006\\
    0.0006    & 0.0004    & 0.0006    & 0.0031
\end{bmatrix}\!,\\[1.5\jot]
C_k &= \setlength\arraycolsep{3pt} \begin{bmatrix}
     1 & 0 & 0 & 0\\
     0 & 0 & 1 & 0
\end{bmatrix}\!,
V_k = \setlength\arraycolsep{3pt} \begin{bmatrix}
     0.0020  & 0.0000\\
     0.0000   & 0.0010
\end{bmatrix}\!,
\end{align*}
for $k \in \mathbb{N}_{[0,N]}$, mean and covariance of the initial condition $m_0 = [0 \ 0 \ 0.2 \  0]^T$ and $M_0 = 10W_k$, weighting coefficients and matrices $Q_{N+1}= \diag\{1,1,1000,1\}$, $\ell_k = 1$, $Q_k= \diag\{1,1,1000,1\}$, and $R_k=1$ for $k \in \mathbb{N}_{[0,N]}$, time horizon $N= 500$, and tradeoff multiplier $\lambda = 0.0066$. For a realization of this system, the value of information, transmission decision, and control input trajectories are shown in Fig.~\ref{c1:fig:trajectories2X}, and the position, velocity, pitch angle, and pitch rate trajectories in Fig.~\ref{c1:fig:trajectories1X}. Note that in this experiment, the value of information became nonnegative only $17$ times, which led to the transmission of a data packet from the sensor to the controller at each of those times. The corresponding trajectories under a periodic triggering policy with the same number of transmissions are also illustrated in Fig.~\ref{c1:fig:trajectories2X} and \ref{c1:fig:trajectories1X}. We observe that the system under the triggering policy designed based on the value of information was able to achieve relatively much better regulation quality.

\section{Conclusion}\label{c1:conclusion}
In this article, we introduced the notion of the value of information as an intrinsic property of networked control systems, and established a theoretical framework for its characterization and computation. The results asserted that the value of information systematically measures the semantics of each data packet as the difference between its benefit and its cost, and that a strategy based on the value of information optimally manages the communication between the sensor and the controller by allowing only data packets with nonnegative valuations to be transmitted. Note that the above objectives could not be achieved by means of the traditional information-theoretic metrics or the traditional event-triggered conditions. We suggest that future research should extend the framework developed in this study to other classes of systems.

\bibliography{../../../mybib}
\bibliographystyle{ieeetr}

\end{document}